\documentclass[12pt, reqno]{amsart}
\usepackage{amssymb, latexsym, mathrsfs, pstcol}

\newtheorem{thm}{Theorem}[section]
\newtheorem{corollary}[thm]{Corollary}
\newtheorem{lemma}[thm]{Lemma}
\newtheorem{theorem}[thm]{Theorem}

\newtheorem{prop}[thm]{Proposition}
\numberwithin{equation}{section}

\theoremstyle{definition}
\newtheorem{rem}[thm]{Remark}

\newcommand{\al}{\alpha}
\renewcommand{\b}{\beta}

\newcommand{\de}{\delta}
\newcommand{\e}{\varepsilon}
\newcommand{\la}{\lambda}
\renewcommand{\phi}{\varphi}

\renewcommand{\d}{\partial}
\newcommand{\R}{{\mathbb R}}

\newcommand{\br}[1]{\left\langle #1 \right\rangle}
\newcommand{\Case}[1]{\noindent \underline{Case #1:}}

\newcommand{\Subcase}[1]{\noindent \underline{Subcase #1:}}

\renewcommand{\qed}{\rule{3mm}{3mm}}
\renewenvironment{proof}
    {\vspace{1mm}\noindent\textbf{Proof.}}
    {\hspace*{\fill} $\qed$\vspace{1mm}}
\newenvironment{proof_of}[1]
    {\vspace{1mm}\noindent {\bf Proof of #1.}}
    {\hspace*{\fill} $\qed$\vspace{1mm}}
\renewcommand{\qed}{\rule{3mm}{3mm}}

\begin{document}
\title[Scattering for Nonlinear Waves with Potential]{Small-Data Scattering for Nonlinear Waves with Potential and
Initial Data of Critical Decay}
\author{Paschalis Karageorgis}
\address{2033 Sheridan Rd., Evanston, IL 60208}
\email{pete@math.northwestern.edu}



\begin{abstract}
We study the scattering problem for the nonlinear wave equation with potential.  In the absence of the potential, one has
sharp global existence results for the Cauchy problem with small initial data; those require the data to decay at a rate
$k\geq k_c$, where $k_c$ is a critical decay rate that depends on the order of the nonlinearity.  However, scattering
results have appeared only for the supercritical case $k>k_c$. In this paper, we extend the latter results to the
critical case and we also allow the presence of a short-range potential.
\end{abstract}
\maketitle

\section{Introduction}
We study the scattering problem for the nonlinear wave equation with potential
\begin{equation}\label{we}
\d_t^2 u - \Delta u + V(x)\cdot u = F(u) \quad\quad \text{in\: $\R^n \times \R$,}
\end{equation}
where $F(u)$ behaves like $|u|^p$ for some $p>1$.  When it comes to the special case $V(x)\equiv 0$, it is known that
both the size of $p$ and the decay rate $k$ of the initial data play a crucial role in the existence theory of the
associated Cauchy problem for small initial data.  In fact, the condition $k\geq 2/(p-1)$ is one of the sharp conditions
needed to ensure the global existence of small-amplitude solutions.  The scattering operator for \eqref{we}, on the other
hand, has been constructed only in the supercritical case $k>2/(p-1)$, where a zero potential was assumed.  In this
paper, we construct the scattering operator for the critical case as well, and we also allow the presence of a small,
rapidly decaying potential.

First, consider the Cauchy problem for \eqref{we} when $V(x)\equiv 0$ and the small initial data are prescribed at time
$t=0$. In what follows, we denote by $p_n$ the positive root of the equation
\begin{equation}\label{pn}
(n-1) p^2 = (n+1)p+2,
\end{equation}
using the convention that $p_1= \infty$.  If $1<p\leq p_n$ and the data are of compact support, then blow-up is known to
generally occur.  And even if $p>p_n$, blow-up may still occur for slowly decaying data of noncompact support. Namely,
under the assumption that
\begin{equation*}
u(x,0)\equiv 0,\quad\quad \d_t u(x,0) \geq \e (1+|x|)^{-k-1},
\end{equation*}
blow-up does occur for any $\e>0$, provided that $0\leq k<2/(p-1)$.  As it turns out, however, the conditions $p>p_n$ and
$k\geq 2/(p-1)$ are not only necessary but also sufficient for the existence of global solutions in the following sense.
Under the assumption that
\begin{equation*}
\d_x^\al u(x,0), \: \d_x^\b \d_t u(x,0)= O(|x|^{-k-1}) \quad\text{as\, $|x|\to \infty$;} \quad |\al|\leq 3,\: |\b|\leq 2,
\end{equation*}
small-amplitude solutions exist for all times when $p>p_n$ and $k\geq 2/(p-1)$, provided that either $n=2,3$; or $n\geq
4$ and the data are of compact support; or $n\geq 4$ and the data are radially symmetric with noncompact support. For a
precise list of references on the Cauchy problem, we refer the reader to \cite{Ka2}.

Next, we focus on the scattering problem for \eqref{we}.  When it comes to the case $V(x)\equiv 0$, the existence of the
scattering operator for small data was established by Pecher \cite{Pe} in $n=3$ space dimensions under the almost optimal
assumptions $p>p_n$ and $k>2/(p-1)$.  Under the exact same assumptions, this result was extended to the case $n=2$ by
Tsutaya \cite{Ts4} and independently by Kubota and Mochizuki \cite{KM}.  As for the higher-dimensional case $n\geq 4$,
the corresponding results were obtained by Kubo and Kubota \cite{KKo, KKe} for radially symmetric data.

In the remaining of this paper, we shall mostly focus on the radially symmetric version of the nonlinear wave equation
with potential \eqref{we}.  Thus, the equation of interest is
\begin{equation}\label{nl}
\left\{
\begin{array}{rll}
\d_t^2 u -\d_r^2 u -\dfrac{n-1}{r} \cdot \d_r u = F(u) - V(r)\cdot u \quad\quad
&\text{in\, $\Omega= (0,\infty)\times \R$}\\
u(r,0) = \phi(r); \quad \d_t u(r,0)= \psi(r) \quad\quad &\text{in\, $\R_+= (0,\infty)$.}
\end{array}\right.
\end{equation}
Before we state our main results, however, let us first introduce some hypotheses.  When it comes to the nonlinear term
$F(u)$, we shall impose the conditions
\begin{equation}\label{F}
F\in \mathcal{C}^1(\R); \quad F(0)= F'(0) = 0; \quad |F'(u)-F'(v)| \leq Ap |u-v|^{p-1}
\end{equation}
for some $A>0$ and some $p>1$.  When it comes to the potential term $V(r)$, we require that
\begin{equation}\label{V}
\sum_{i=0}^1 \br{r}^i \,|V^{(i)}(r)| \leq V_0 \br{r}^{-\kappa}, \quad\quad \kappa>2
\end{equation}
for some small $V_0>0$, where the bracket notation $\br{s}= 1+|s|$ is used for convenience.  As for the initial data, we
shall assume that
\begin{equation}\label{da3}
\sum_{i=0}^2 \br{r}^i \,|\phi^{(i)}(r)| + \sum_{i=0}^1 \br{r}^{i+1} \,|\psi^{(i)}(r)| \leq \e\br{r}^{-k}
\end{equation}
for some small $\e>0$ and some $k> 0$.

Although our approach does not really depend on the parity of $n$, the decay estimates that one can obtain for solutions
to \eqref{nl} inevitably do.  In particular, it will be convenient to introduce the parameters
\begin{equation}\label{am}
(a,m) = \left\{
\begin{array}{ccl}
\left( 1\:,\: \frac{n-3}{2} \right) &&\text{if\, $n$ is odd} \\
\left( \frac{1}{2} \:,\: \frac{n-2}{2} \right) &&\text{if\, $n$ is even}
\end{array}\right.
\end{equation}
we shall frequently use in what follows.  Let us remark that $m\geq 1$ when $n\geq 4$ and that the sum $a+m= (n-1)/2$ is
independent of the parity of $n$.

Our plan is to construct a solution of \eqref{nl} that is continuously differentiable and belongs to the Banach space
\begin{equation}\label{X}
X= \left\{ u(r,t)\in \mathcal{C}^1(\Omega) \::\: ||u|| < \infty \right\}, \quad\quad \Omega= (0,\infty) \times \R.
\end{equation}
Here, the norm $||\cdot||$ is defined by
\begin{equation}\label{N}
||u|| = \sum_{j=0}^1 \sup_{(r,t)\in \Omega} \: |\d_r^j u(r,t)| \cdot r^{m-1+j} \br{r}^{1-j} \cdot W_k(r,|t|),
\end{equation}
where $m$ is given by \eqref{am} and the weight function $W_k$ is of the form
\begin{equation}\label{W1}
W_k(r,|t|) = \br{|t|+r}^a \br{|t|-r}^\nu,\quad\quad \nu=k-m-a.
\end{equation}
This weighted norm is partly dictated by the available estimates regarding the solutions to the homogeneous problem
\begin{equation}\label{he}
\d_t^2 u_0 - \d_r^2 u_0  - \dfrac{n-1}{r} \cdot \d_r u_0 = 0 \quad\quad\text{in\, $\Omega= \R_+\times \R$.}
\end{equation}
For a proof of the following lemma, we refer the reader to \cite{Ka1, KKo, KKe}.

\begin{lemma}\label{dec}
Let $n\geq 4$ be an integer and define $a,m$ by \eqref{am}.  Suppose that $\phi\in \mathcal{C}^2(\R_+)$ and $\psi\in
\mathcal{C}^1(\R_+)$ are subject to \eqref{da3} for some $\e>0$ and some $(n-1)/2< k <n-1$. Then the homogeneous equation
\eqref{he} subject to initial data $u_0(r,0)=\phi(r)$, $\d_t u_0(r,0)= \psi(r)$ admits a unique solution $u_0^-\in
\mathcal{C}^1(\Omega)$ which satisfies
\begin{equation}\label{en}
|\d_r^j u_0^-(r,t)| + |\d_t^j u_0^-(r,t)| \leq C_0(k,n) \cdot \e r^{1-m-j} \br{r}^{j-1} \cdot W_k(r,|t|)^{-1}
\end{equation}
when $j=0,1$.  In particular, $u_0^-$ is in the Banach space \eqref{X} and we have $||u_0^-||\leq C_0\e$.
\end{lemma}

\begin{rem}
An analogous result holds when $n=2,3$, in which case the estimate
\begin{equation}\label{N3}
\sum_{|\b|\leq 2} |\d_x^\b u_0^-(x,t)| \leq C(k) \cdot \e W_k(|x|,|t|)^{-1}
\end{equation}
holds even without the assumption of radial symmetry; see \cite{Ts1,Ts3}, for instance.
\end{rem}

Regarding the solutions to the nonlinear equation, the existence result of this paper can now be stated as follows.

\begin{theorem}\label{et}
Let $n\geq 4$ be an integer.  Let $\phi\in \mathcal{C}^2(\R_+)$ and $\psi\in \mathcal{C}^1(\R_+)$ satisfy \eqref{da3} for
some $\e,k>0$. Suppose the nonlinear term $F(u)$ satisfies \eqref{F} for some $A>0$ and some
\begin{equation}\label{pf}
p_n < p < 1 + \frac{4}{n-1} \:,
\end{equation}
where $p_n$ is the positive root of \eqref{pn}.  Also, assume $V(r)$ is subject to \eqref{V} for some $V_0>0$ and
$\kappa>2$. If $k\geq 2/(p-1)$ and if $V_0,\e$ are sufficiently small, then the following conclusions hold with $D=
(\d_r, \d_t)$ and $u_0^-$ as in Lemma \ref{dec}.

\begin{itemize}
\item[(a)]
There exists a unique solution $u\in X$ to the nonlinear equation \eqref{nl} such that
\begin{equation*}
|D^\b (u(r,t)-u_0^-(r,t))| \leq C||u|| \cdot r^{1-|\b|-m} \br{r}^{|\b|-1} \cdot W_k(r,|t|)^{-1}, \quad\quad |\b|\leq 1
\end{equation*}
for some constant $C$ which is independent of $r,t$.  Besides, one has $||u||\leq 1$.

\item[(b)]
There exists a unique solution $u_0^+\in X$ to the homogeneous equation \eqref{he} such that
\begin{equation*}
|D^\b (u(r,t)-u_0^+(r,t))| \leq C||u|| \cdot r^{1-|\b|-m} \br{r}^{|\b|-1} \cdot W_k(r,|t|)^{-1}, \quad\quad |\b|\leq 1
\end{equation*}
for some constant $C$ which is independent of $r,t$.
\end{itemize}
\end{theorem}

\begin{rem}
Using the norm dictated by \eqref{N3} instead of our norm \eqref{N}, one can extend Theorem \ref{et} to the case $n=3$
without imposing radial symmetry.  When $n=1,2$, on the other hand, there exist arbitrarily small, rapidly decaying
potentials for which small-data blow-up occurs; see section 6 in \cite{Ka2} for more details. Nevertheless, it is still
possible to treat the case $n=2$ when $V(x)\equiv 0$; for supercritical decay rates, we refer the reader to \cite{Ts4}.
\end{rem}

\begin{rem}
Under our assumption that $k\geq 2/(p-1)$, the upper bound in \eqref{pf} implies
\begin{equation}\label{kl}
k\geq \frac{2}{p-1} > \frac{n-1}{2} =a+m,
\end{equation}
where $a,m$ are defined by \eqref{am}.  It is possible to weaken the upper bound \eqref{pf} on $p$ in the spirit of
\cite{Ka2} and thus obtain existence for decay rates $k\leq (n-1)/2$ as well.  For such decay rates, however, we are
unable to show that the solutions are asymptotically free, as we do in our next corollary.
\end{rem}

\begin{corollary}\label{sc}
Let the assumptions of the previous theorem hold and define the energy norm
\begin{equation}\label{ener}
||w(\cdot\,,t)||_e = \left( \int_0^\infty \left( |\d_r w(r,t)|^2 + |\d_t w(r,t)|^2 \right) \cdot r^{n-1} \:dr
\right)^{1/2}.
\end{equation}
If it happens that $k> m+1$, then one has
\begin{align}
||u(\cdot\,,t) - u_0^-(\cdot\,,t)||_e &\leq C||u||\cdot \br{t}^{-\theta}, \quad\quad t\leq 0; \label{sc1}\\
||u(\cdot\,,t) - u_0^+(\cdot\,,t)||_e &\leq C||u||\cdot \br{t}^{-\theta}, \quad\quad t\geq 0.\label{sc2}
\end{align}
Here, the constant $C$ is independent of $t$ and we have also set
\begin{equation}\label{th}
\theta \equiv \min \bigl( (a+m)(p-1)-1,\: (k-m)p-1,\: k-m -1\bigr) >0
\end{equation}
for convenience.
\end{corollary}

The remaining of this paper is organized as follows.  In section \ref{pre}, we collect some facts about the inhomogeneous
wave equation and reduce the proof of our existence result to the estimation of certain integrals involving our weight
function \eqref{W1}.  In section \ref{slow}, we gather the necessary estimates to treat those integrals, while section
\ref{ex} is devoted to the proofs of our main results, Theorem \ref{et} and Corollary \ref{sc}.

\section{Preliminaries}\label{pre}
In this section, we prepare a few basic lemmas that will be needed in the proof of our existence theorem regarding the
nonlinear wave equation with potential
\begin{equation}\label{nlf}
\d_t^2 u -\d_r^2 u -\dfrac{n-1}{r} \cdot \d_r u = F(u) - V(r)\cdot u \quad\quad\text{in\, $\Omega= (0,\infty)\times \R$.}
\end{equation}
Recall that we seek a solution to \eqref{nlf} for initial data of decay rate $k\geq 2/(p-1)$. Since there is no loss of
generality in decreasing this decay rate, we may take $k$ to be smaller than any quantity that exceeds $2/(p-1)$. Now,
our assumption \eqref{pf} that $p$ is larger than the critical power $p_n$ of \eqref{pn} can be written as
\begin{equation*}
\frac{2}{p-1} < \frac{n-1}{2} \cdot p - 1 = (a+m)p -1
\end{equation*}
with $a,m$ as in \eqref{am}; and it can also be written as
\begin{equation*}
\frac{2}{p-1} < \frac{n-1}{2} + \frac{1}{p} = a+m + \frac{1}{p} \:.
\end{equation*}
This allows us to decrease the value of $k$ and henceforth assume
\begin{equation}\label{k1}
\frac{2}{p-1} \leq k < \min((a+m)p -1, \:a+m+1/p)
\end{equation}
without loss of generality.  Since we also have $p>1$, this automatically implies
\begin{equation}\label{k2}
k < a+m+1 = \frac{n+1}{2} \:.
\end{equation}
Finally, it is convenient to decrease the decay rate $\kappa >2$ of the potential $V(r)$ so that
\begin{equation}\label{ka}
\kappa < m+2.
\end{equation}
We can do this without loss of generality whenever $m>0$, namely, whenever $n\geq 4$.

\begin{rem}
The last simplification is not available when $n=2,3$, yet we shall only need it to control certain integrals that do not
arise in odd dimensions due to the strong form of Huygens' principle; see Lemma \ref{J2}.  This is exactly the point
where our proof would fail to work in two space dimensions, yet this is merely because of the potential term.
\end{rem}

The following elementary fact is quite standard.

\begin{lemma}\label{new}
Let $L$ denote the Riemann operator for the wave equation in the radial case.  Given a function $G$ of two variables, we
define the Duhamel operator $\mathscr{L}$ as
\begin{equation}\label{Du}
[\mathscr{L}G](r,t) = \int_{-\infty}^t [LG(\cdot\,,\tau)](r,t-\tau) \:d\tau.
\end{equation}
When $G\in \mathcal{C}^1(\Omega)$, one then has $\mathscr{L}G\in \mathcal{C}^1(\Omega)$ and this function provides a
solution to
\begin{equation*}
\left( \d_t^2 - \d_r^2 - \frac{n-1}{r} \cdot \d_r \right) [\mathscr{L}G](r,t) = G(r,t) \quad\quad \text{in\, $\Omega=
(0,\infty) \times \R$.}
\end{equation*}
\end{lemma}

Although we shall not need an explicit representation for the Riemann operator here, the following basic estimate for the
Duhamel operator is important.  Its proof essentially repeats that of Proposition 2.3 in \cite{Ka2}, so we are going to
omit it.

\begin{prop}\label{Les}
Let $n\geq 4$ be an integer and define $a,m$ by \eqref{am}.  Suppose $G\in \mathcal{C}^1(\Omega)$ satisfies the
singularity condition
\begin{equation}\label{sg2}
G(\la,\tau) = O \left( \la^{-2m-2+\de} \right) \quad \text{as $\la \rightarrow 0$}
\end{equation}
for some fixed $\de>0$.  With $D= (\d_r, \d_t)$ and $\la_\pm = t-\tau\pm r$, one then has
\begin{align*}
|D^\b [\mathscr{L}G](r,t)|
&\leq Cr^{j-|\b|-m-a} \int_{-\infty}^t \int_{|\la_-|}^{\la_+} \frac{\la^{m-j+1}}{(\la-\la_-)^{1-a}} \cdot \sum_{s=0}^j
\la^s \,|\d_\la^s G(\la,\tau)| \:d\la\,d\tau \\
&\quad + Cr^{j-|\b|-m} \int_{-\infty}^{t-r} \int_0^{\la_-} \frac{\la^{2m-j+1}}{\la_-^m \la_+^a \,(\la_- -\la)^{1-a}} \cdot
\sum_{s=0}^j \la^s \,|\d_\la^s G(\la,\tau)| \:d\la\,d\tau\\
&\quad + Cr^{j-|\b|-m-a} \int_{t-2r}^t  |\la_\pm|^{a+m-j+1} \left[ \sum_{s=0}^{j-1} \la^s \,|\d_\la^s G(\la,\tau)|
\right]_{\la= |\la_\pm|} d\tau
\end{align*}
whenever $(r,t)\in \Omega$ and $|\b|\leq j\leq 1$.  Besides, the constant $C$ is independent of $r,t$.
\end{prop}

\begin{lemma}\label{bas}
Suppose $u$ belongs to the Banach space \eqref{X} and let $p>1$.  Assuming \eqref{F} and \eqref{V}, one then has
\begin{equation}\label{ma1}
\sum_{s=0}^{j_0} \la^s \:|\d_\la^s F(u(\la,\tau))| \leq 2Ap ||u||^p \cdot\la^{j-mp} \br{\la}^{j_0-j} \cdot
W_k(\la,|\tau|)^{-p}
\end{equation}
and also
\begin{equation}\label{ma2}
\sum_{s=0}^j \la^s \:|\d_\la^s (V(\la)\cdot u(\la,\tau))| \leq 4V_0 ||u|| \cdot \la^{j-1-m} \br{\la}^{1-\kappa} \cdot
W_k(\la,|\tau|)^{-1}
\end{equation}
whenever $(\la,\tau)\in \Omega$ and $0\leq j, j_0\leq 1$.
\end{lemma}

\begin{proof}
Because of our assumption \eqref{F}, the fundamental theorem of calculus ensures that
\begin{equation*}
|F(u)| \leq A\cdot |u|^p, \quad\quad |F'(u)|\leq Ap\cdot |u|^{p-1}.
\end{equation*}
In particular, it ensures that
\begin{equation*}
|\d_\la^s F(u(\la,\tau))| \leq Ap^s \cdot |u(\la,\tau)|^{p-s} \cdot |\d_\la u(\la,\tau)|^s, \quad\quad s=0,1.
\end{equation*}
Recalling the definition \eqref{N} of our norm, the last equation easily leads to
\begin{equation*}
\la^s \,|\d_\la^s F(u(\la,\tau))|\leq Ap^s ||u||^p \cdot\la^{p-mp}\br{\la}^{s-p}\cdot W_k(\la,|\tau|)^{-p}, \quad\quad
s=0,1.
\end{equation*}
In view of our assumptions that $j\leq 1< p$, this also implies
\begin{equation*}
\sum_{s=0}^{j_0} \la^s \,|\d_\la^s F(u(\la,\tau))| \leq Ap ||u||^p \cdot \sum_{s=0}^{j_0} \la^{j-mp} \br{\la}^{s-j} \cdot
W_k(\la,|\tau|)^{-p}, \quad\quad j_0=0,1.
\end{equation*}
Moreover, $s-j\leq j_0-j$ within the last sum, so our first assertion \eqref{ma1} follows.

Since our second assertion \eqref{ma2} is easier to establish, we shall omit the details.
\end{proof}

\section{A Priori Estimates}\label{slow}
In this section, we prepare the main estimates needed in the proof of our existence result.  Those are contained in the
following

\begin{theorem}\label{bes}
Let $n\geq 4$ be an integer and $a,m$ be as in \eqref{am}.  Suppose $F(u)$ satisfies \eqref{F} for some $p$ subject to
\eqref{pf} and that $V(r)$ satisfies \eqref{V}.  Assume the decay rates $k,\kappa$ are subject to \eqref{kl}, \eqref{k1}
and \eqref{ka}.  Define the Duhamel operator $\mathscr{L}$ by \eqref{Du} and let $u$ be an element of the Banach space
\eqref{X}. With $D= (\d_r, \d_t)$, one then has
\begin{equation}\label{mb1}
|D^\b [\mathscr{L} F(u)](r,t)| \leq C_1 ||u||^p \cdot r^{1-|\b|-m} \br{r}^{|\b|-1} \cdot W_k(r,|t|)^{-1}
\end{equation}
and also
\begin{equation}\label{mb2}
|D^\b [\mathscr{L} (Vu)](r,t)| \leq C_1 V_0 ||u||\cdot r^{1-|\b|-m} \br{r}^{|\b|-1} \cdot W_k(r,|t|)^{-1}
\end{equation}
as long as $(r,t)\in \Omega$ and $|\b|\leq 1$.  Besides, the constant $C_1$ is independent of $r, t$.
\end{theorem}

\begin{rem}
As an immediate consequence of \eqref{pf}, one finds that the conditions
\begin{equation}\label{con}
0< a(p-1) < 2 - mp+ m < ap
\end{equation}
also hold under the assumptions of this theorem.
\end{rem}

The proof of Theorem \ref{bes} is given in the next section.  Here, we shall merely study certain integrals which arise
in the course of the proof. Those involve our weight function \eqref{W1} and some other parameters we have introduced
\eqref{am}. Throughout this section, in particular, we shall assume
\begin{equation}\label{ks}
0< \nu < 1/p,\quad \nu= k-m-a; \quad m\geq 0; \quad a>0; \quad \kappa>2;\quad p>1.
\end{equation}
Let us remark that the two leftmost inequalities follow by \eqref{kl} and \eqref{k1}, respectively.

\begin{lemma}\label{ts1}
Let $r,t>0$ be arbitrary.  Assuming that $a,\nu>0$, one has
\begin{equation*}
A_1 \equiv \int_{|t-r|}^{t+r} \frac{\br{y}^{-a-\nu} \,dy}{(r-t+y)^{1-a}} \leq C(a,\nu) \cdot r^a W_k(r,t)^{-1},
\end{equation*}
where the weight function $W_k$ is given by \eqref{W1}.
\end{lemma}

\begin{proof}
If either $r\leq 1$ or $t\geq 2r$, then $\br{t-r}$ is equivalent to $\br{t+r}$, so we get
\begin{equation*}
A_1 \leq C\br{t\pm r}^{-a-\nu} \int_{|t-r|}^{t+r} (r-t+y)^{a-1} \:dy \leq Cr^a \cdot \br{t\pm r}^{-a-\nu}
\end{equation*}
because $a>0$.  This does imply the desired estimate under the present assumptions.

If $r\geq 1$ and $t\leq 2r$, on the other hand, an integration by parts gives
\begin{equation*}
A_1 = \left[ \frac{1}{a} \: (r-t+y)^a \cdot \br{y}^{-a-\nu} \right]_{y=|t-r|}^{t+r} + \frac{a+\nu}{a} \int_{|t-r|}^{t+r}
(r-t+y)^a \cdot \br{y}^{-a-\nu-1} \:dy.
\end{equation*}
Since $a>0$ and since $r-t+y \leq 2y$ within the region of integration, we now get
\begin{equation*}
A_1 \leq Cr^a \br{t+r}^{-a-\nu} + C\int_{|t-r|}^{t+r} \br{y}^{-\nu-1} \:dy \leq C\br{t-r}^{-\nu}
\end{equation*}
because $\nu>0$.  As $r$ is equivalent to $\br{t+r}$ when $r\geq \max(t/2,1)$, the result follows.
\end{proof}

\begin{lemma}\label{fa1}
Let $z\geq 0$ be arbitrary.  Given constants $a>0$ and $b<1$, one then has
\begin{equation}\label{fa1a}
A_2 \equiv \int_0^z \frac{\br{y}^{-b} \,dy}{(z\pm y)^{1-a}} \leq C(a,b) \cdot z^a \br{z}^{-b}.
\end{equation}
\end{lemma}

\begin{proof}
If it happens that $0\leq z\leq 1$, then we easily get
\begin{equation*}
A_2 \leq C\int_0^z (z\pm y)^{a-1}\:dy = C'z^a
\end{equation*}
because $a>0$.  Moreover, this does imply the desired \eqref{fa1a} whenever $z$ is bounded.

If it happens that $z\geq 1$, on the other hand, we get
\begin{equation*}
A_2 \leq Cz^{a-1} \int_0^{z/2} \br{y}^{-b} \,dy + C\br{z}^{-b} \int_{z/2}^z (z\pm y)^{a-1} \:dy.
\end{equation*}
In view of our assumptions that $a>0$ and $b<1$, the desired estimate \eqref{fa1a} follows.
\end{proof}

\begin{lemma}\label{fa2}
Let $y\in \R$ and $z\geq |y|$.  Assuming \eqref{k1}, \eqref{con} and \eqref{ks}, one has
\begin{equation*}
I \equiv \int_z^\infty \br{x+y}^{1-\kappa} \cdot \br{x}^{-a} \,dx \leq C(\kappa)\cdot \br{z}^{-a}
\end{equation*}
and also
\begin{equation}\label{fa2a}
J \equiv \int_z^\infty (x+y)^{1-mp+m} \cdot \br{x}^{-ap} \,dx \leq C(a,m,p) \cdot \br{z}^{\nu p -\nu-a}.
\end{equation}
\end{lemma}

\begin{proof}
We only prove our second assertion, as our first assertion is trivial.  Write
\begin{align*}
J &= \int_z^{2z+1} (x+y)^{1-mp+m} \cdot \br{x}^{-ap} \,dx + \int_{2z+1}^\infty (x+y)^{1-mp+m} \cdot \br{x}^{-ap} \,dx\\
&\equiv J_1 + J_2.
\end{align*}
When it comes to the first part, $\br{x}$ is equivalent to $\br{z}$ and we easily get
\begin{equation*}
J_1 \leq C\br{z}^{-ap} \int_z^{2z+1} (x+y)^{1-mp+m} \,dx \leq C\br{z}^{2-mp+m-ap}
\end{equation*}
because $2-mp+m>0$ by \eqref{con}.  Since we also have
\begin{equation}\label{ide}
2- mp +m - ap = [2-k(p-1)] + \nu (p-1) -a \leq \nu p -\nu -a
\end{equation}
by \eqref{ks} and \eqref{k1}, we may combine the last two equations to deduce the desired \eqref{fa2a}.

When it comes to $J_2$, on the other hand, we have $2x\geq x+ 2|y|\geq x+y\geq 1$ within the region of integration, so
\begin{equation*}
J_2 \leq C\int_{2z+1}^\infty \br{x+y}^{1-mp+m-ap} \:dx \leq C\br{z}^{2-mp+m-ap}
\end{equation*}
because $2-mp+m-ap<0$ by \eqref{con}.  In view of \eqref{ide}, the desired \eqref{fa2a} now follows.
\end{proof}

\begin{lemma}\label{I1}
Let $(r,t)\in \Omega$.  Assuming \eqref{k1}, \eqref{con} and \eqref{ks}, one has
\begin{equation*}
\mathcal{I}_1 \equiv \int_{-\infty}^t \int_{|\la_-|}^{\la_+} \frac{\br{\la}^{1-\kappa} \cdot
W_k(\la,|\tau|)^{-1}}{(\la-\la_-)^{1-a}} \:\:d\la\,d\tau \leq Cr^a W_k(r,|t|)^{-1},
\end{equation*}
where $\la_\pm = t-\tau\pm r$, $W_k$ is given by \eqref{W1} and the constant $C$ is independent of $r, t$.
\end{lemma}

\begin{proof}
According to Lemma 3.5 in \cite{Ka2}, we do have the estimate
\begin{equation*}
\int_{\min (0,t)}^t \int_{|\la_-|}^{\la_+} \frac{\br{\la}^{1-\kappa} \cdot W_k(\la,|\tau|)^{-1}}{(\la-\la_-)^{1-a}}
\:\:d\la\,d\tau \leq Cr^a W_k(r,|t|)^{-1}
\end{equation*}
when $t\geq 0$.  Since this estimate holds trivially when $t\leq 0$, it remains to show that
\begin{equation}\label{sh1}
\mathcal{I}_1' \equiv \int_{-\infty}^{\min (0,t)} \int_{|\la_-|}^{\la_+} \frac{\br{\la}^{1-\kappa} \cdot
W_k(\la,|\tau|)^{-1}}{(\la-\la_-)^{1-a}} \:\:d\la\,d\tau \leq Cr^a W_k(r,|t|)^{-1}.
\end{equation}
Let us first recall our definition \eqref{W1} of our weight function $W_k$ and write
\begin{equation*}
\mathcal{I}_1' = \int_{-\infty}^{\min (0,t)} \int_{|\la_-|}^{\la_+} \frac{\br{\la}^{1-\kappa} \cdot \br{\la-\tau}^{-a}
\br{\la+\tau}^{-\nu}}{(r-t+\la+\tau)^{1-a}} \:\:d\la\,d\tau.
\end{equation*}
As $\tau\leq t$ within the region of integration, we have $\la\geq |\la_-|\geq t-\tau -r$.  As $\tau\leq 0$, we also have
$\la\geq |\la_-|\geq |t-r| + \tau$.  Changing variables by $x= \la-\tau$ and $y=\la+\tau$, we then get
\begin{align*}
\mathcal{I}_1' &\leq C\int_{t-r}^{t+r} \frac{\br{y}^{-\nu}}{(r-t+y)^{1-a}} \int_{\max(|y|,|t-r|)}^\infty
\br{x+y}^{1-\kappa} \cdot \br{x}^{-a} \:dx\,dy \\
&\leq C\int_{t-r}^{t+r} \frac{\br{y}^{-\nu}}{(r-t+y)^{1-a}} \cdot \max(\br{y}, \br{t-r})^{-a} \:dy
\end{align*}
by Lemma \ref{fa2} with $z= \max(|y|, |t-r|)$.

Note that $\nu< \nu p <1$ by \eqref{ks}.  To show that the last integral satisfies the desired \eqref{sh1}, we shall now
establish the more general estimate
\begin{equation}\label{mg}
\mathcal{I}_1'' \equiv \int_{t-r}^{t+r} \frac{\br{y}^{-b_1}}{(r-t+y)^{1-a}} \cdot \max(\br{y}, \br{t-r})^{-b_2} \:dy \leq
Cr^a W_k(r,|t|)^{-1}
\end{equation}
for any constants $b_1,b_2\in \R$ with $b_1<1$ and $b_1+b_2= a+\nu$.

\Case{1} When $t\geq 0$, we have $t-r\leq |t-r|\leq t+r$, so our definition \eqref{mg} reads
\begin{align*}
\mathcal{I}_1'' &= \br{t-r}^{-b_2} \int_{t-r}^{|t-r|} \frac{\br{y}^{-b_1}}{(r-t+y)^{1-a}} \:\:dy + \int_{|t-r|}^{t+r}
\frac{\br{y}^{-a-\nu}}{(r-t+y)^{1-a}} \:\:dy \\
&\equiv \mathcal{I}_{11}'' + \mathcal{I}_{12}''.
\end{align*}
To treat the latter integral, we need only invoke Lemma \ref{ts1}.  To treat the former integral, we may assume that
$r\geq t$.  In view of our assumption that $b_1<1$, we then get
\begin{equation*}
\mathcal{I}_{11}'' \leq 2\br{r-t}^{-b_2} \int_0^{r-t} \frac{\br{y}^{-b_1}}{(r-t\pm y)^{1-a}} \:\:dy \leq C(r-t)^a \cdot
\br{r-t}^{-a-\nu}
\end{equation*}
by Lemma \ref{fa1} with $z=r-t\geq 0$.  In particular, it remains to show that
\begin{equation*}
(r-t)^a \cdot \br{r-t}^{-a-\nu} \leq Cr^a \cdot \br{r+t}^{-a} \br{r-t}^{-\nu} \quad\quad\text{whenever \,$r\geq t\geq
0$.}
\end{equation*}
If $r\geq t$ and $r\leq 1$, this is easy to see because $r-t\leq r$ and $\br{r-t}$ is equivalent to $\br{r+t}$.  If
$r\geq t$ and $r\geq 1$, on the other hand, $r$ is equivalent to $\br{r+t}$ and we similarly get
\begin{equation*}
(r-t)^a \cdot \br{r-t}^{-a-\nu} \leq \br{r-t}^{-\nu} \leq Cr^a \cdot \br{r+t}^{-a} \br{r-t}^{-\nu}.
\end{equation*}

\Case{2} When $t\leq 0$, we have $|r+t|\leq r+|t|=r-t$, so our definition \eqref{mg} reads
\begin{equation}\label{le}
\mathcal{I}_1'' = \br{r-t}^{-b_2}  \int_{t-r}^{t+r} \frac{\br{y}^{-b_1}}{(r-t+y)^{1-a}} \:\:dy.
\end{equation}

\Subcase{2a} If it happens that $|t|\leq 3r$, we proceed as in the previous case to obtain
\begin{equation*}
\mathcal{I}_1'' \leq \br{r-t}^{-b_2} \int_{t-r}^{r-t} \frac{\br{y}^{-b_1}}{(r-t+y)^{1-a}} \:\:dy \leq C(r-t)^a \cdot
\br{r-t}^{-a-\nu}
\end{equation*}
using Lemma \ref{fa1}.  Since $r-t= r+|t|\leq 4r$ for this subcase, the result follows easily.

\Subcase{2b} If it happens that $-t=|t|\geq 3r$, then $\br{r+t}$ is equivalent to $\br{r-t}$ because
\begin{equation*}
|r+t|\leq r+|t| = r-t \leq -2(r+t)
\end{equation*}
for this subcase.  In particular, equation \eqref{le} trivially leads to
\begin{equation*}
\mathcal{I}_1'' \leq C\br{r\pm t}^{-a-\nu} \int_{t-r}^{t+r} (r-t+y)^{a-1} \:dy \leq Cr^a \br{r\pm t}^{-a-\nu}
\end{equation*}
since $a>0$.  This does imply the desired \eqref{mg} whenever $\br{r+t}$ is equivalent to $\br{r-t}$.
\end{proof}

\begin{lemma}\label{J1}
Let $(r,t)\in \Omega$.  Assuming \eqref{k1}, \eqref{con} and \eqref{ks}, one has
\begin{equation*}
\mathcal{J}_1 \equiv \int_{-\infty}^t \int_{|\la_-|}^{\la_+} \frac{\la^{1-mp+m} \cdot W_k(\la,|\tau|)^{-p}}{(\la
-\la_-)^{1-a}} \:\:d\la\,d\tau \leq Cr^a W_k(r,|t|)^{-1},
\end{equation*}
where $\la_\pm = t-\tau\pm r$, $W_k$ is given by \eqref{W1} and the constant $C$ is independent of $r, t$.
\end{lemma}

\begin{proof}
According to Lemma 3.6 in \cite{Ka2}, we do have the estimate
\begin{equation*}
\int_{\min (0,t)}^t \int_{|\la_-|}^{\la_+} \frac{\la^{1-mp+m} \cdot W_k(\la,|\tau|)^{-p}}{(\la -\la_-)^{1-a}}
\:\:d\la\,d\tau \leq Cr^a W_k(r,|t|)^{-1}
\end{equation*}
when $t\geq 0$.  Since this estimate holds trivially when $t\leq 0$, it remains to show that
\begin{equation*}
\mathcal{J}_1' \equiv \int_{-\infty}^{\min (0,t)} \int_{|\la_-|}^{\la_+} \frac{\la^{1-mp+m} \cdot
W_k(\la,|\tau|)^{-p}}{(\la -\la_-)^{1-a}} \:\:d\la\,d\tau \leq Cr^a W_k(r,|t|)^{-1}.
\end{equation*}
Let us now proceed as in the proof of the previous lemma. Changing variables by $x= \la-\tau$ and $y=\la+\tau$, we use
Lemma \ref{fa2} to arrive at
\begin{align*}
\mathcal{J}_1' &\leq C\int_{t-r}^{t+r} \frac{\br{y}^{-\nu p}}{(r-t+y)^{1-a}} \int_{\max(|y|,|t-r|)}^\infty
(x+y)^{1-mp+m} \cdot \br{x}^{-ap} \:dx\,dy \\
&\leq C\int_{t-r}^{t+r} \frac{\br{y}^{-\nu p}}{(r-t+y)^{1-a}} \cdot \max(\br{y}, \br{t-r})^{\nu p-\nu-a} \:dy.
\end{align*}
Since $\nu p<1$ by \eqref{ks}, we may then invoke our estimate \eqref{mg} to complete the proof.
\end{proof}

\begin{lemma}\label{fa5}
Let $y\in \R$.  Assuming that $\kappa>2$ and $0\leq b<1$, one has
\begin{align*}
A_3 &\equiv \int_{-y}^\infty \br{x+y}^{1-\kappa} \cdot \br{x}^{-b} \:dx \leq C(b,\kappa)\cdot \br{y}^{-b}.
\end{align*}
Assuming that $2-mp+m>0$ and that $b<1< m(p-1)+b$, one also has
\begin{align*}
A_4 &\equiv \int_{-y}^\infty (x+y)^{1-mp+m} \cdot \br{x+y}^{-1} \br{x}^{-b} \:dx \leq C(b,m,p)\cdot \br{y}^{2-mp+m-b}.
\end{align*}
\end{lemma}

\begin{proof}
We shall only prove our second assertion, as the proof of our first assertion is quite similar.  When $x\geq 2|y|+1$, one
has $2\leq x+1\leq 2(x+y)\leq 3(x+1)$ and this implies
\begin{equation*}
A_4\leq C\int_{2|y|+1}^\infty \br{x}^{-m(p-1)-b} \:dx + \int_{-y}^{2|y|+1} (x+y)^{1-mp+m} \cdot \br{x+y}^{-1}
\br{x}^{-b} \:dx.
\end{equation*}
Using our assumption that $m(p-1)+b>1$, we now arrive at
\begin{equation*}
A_4\leq C\br{y}^{2-m(p-1)-b} + \int_{-y}^{2|y|+1} (x+y)^{1-mp+m} \cdot \br{x}^{-b} \:dx.
\end{equation*}
Since $\br{x}$ is equivalent to $\br{y}$ whenever $|y|\leq x\leq 2|y|+1$, this easily leads us to
\begin{equation*}
A_4\leq C\br{y}^{2-m(p-1)-b} + \int_{-y}^{|y|} (x+y)^{1-mp+m} \cdot \br{x}^{-b} \:dx
\end{equation*}
because $2-mp+m>0$ by assumption.  To finish the proof, it thus remains to show that
\begin{equation}\label{sh4}
A_4' \equiv \int_{-y}^{y} (x+y)^{1-mp+m} \cdot \br{x}^{-b} \:dx \leq C\br{y}^{2-m(p-1)-b}
\end{equation}
whenever $y\geq 0$.  If it happens that $0\leq y\leq 1$, this estimate is easy to establish.  Let us then assume that
$y\geq 1$.  Since $\br{x}$ is equivalent to $\br{y}$ for each $-y\leq x\leq -y/2$ and since $x+y$ is equivalent to
$\br{y}$ for each $-y/2\leq x\leq y$, we find
\begin{equation*}
A_4'\leq C\br{y}^{-b} \int_{-y}^{-y/2} (x+y)^{1-mp+m} \:dx + C\br{y}^{1-mp+m} \int_{-y/2}^y \br{x}^{-b} \:dx.
\end{equation*}
In view of our assumptions that $2-mp+m>0$ and $b<1$, the desired \eqref{sh4} follows.
\end{proof}

\begin{rem}
In what follows, we shall need to apply Lemma \ref{fa5} when $b=a(p-1)$ and also when $b=\nu p$.  To ensure the lemma is
applicable, we shall thus need to know that
\begin{equation}\label{iwcs}
a(p-1) < 1 < (a+m)(p-1), \quad\quad \nu p< 1< m(p-1) +\nu p.
\end{equation}
These inequalities follow from \eqref{pf}, \eqref{k1} and \eqref{ks} when $a,m$ are defined by \eqref{am} for some
integer $n\geq 4$. For instance, \eqref{k1} and \eqref{ks} combine to give
\begin{equation*}
(a+m)(p-1) = [(a+m)p -k] + [k-m-a] > 1.
\end{equation*}
Since the remaining assertions of \eqref{iwcs} are also easy to verify, we shall omit the details.
\end{rem}

\begin{lemma}\label{Ipm}
Let $(r,t)\in \Omega$.  Assuming \eqref{k1}, \eqref{con} and \eqref{ks}, one has
\begin{equation}\label{Ipm1}
\mathcal{I}_\pm \equiv \int_{t-2r}^t |\la_\pm|^a \cdot \br{\la_\pm}^{1 -\kappa} \cdot W_k(|\la_\pm|, |\tau|)^{-1} \:d\tau
\leq Cr^a W_k(r,|t|)^{-1},
\end{equation}
where $\la_\pm = t-\tau \pm r$, $W_k$ is given by \eqref{W1} and the constant $C$ is independent of $r,t$.
\end{lemma}

\begin{proof}
According to Lemma 3.12 in \cite{Ka2}, we do have the estimate
\begin{equation*}
\int_{\max(t-2r,0)}^t |\la_\pm|^a \cdot \br{\la_\pm}^{1 -\kappa} \cdot W_k(|\la_\pm|, |\tau|)^{-1} \:d\tau \leq Cr^a
W_k(r,|t|)^{-1}
\end{equation*}
whenever $t\geq 0$.  Instead of \eqref{Ipm1}, it thus suffices to show that
\begin{equation}\label{Ipm2}
\mathcal{I}_\pm' \equiv \int_{t-2r}^{\min(t,0)} |\la_\pm|^a \cdot \br{\la_\pm}^{1 -\kappa} \cdot W_k(|\la_\pm|,
|\tau|)^{-1} \:d\tau \leq Cr^a W_k(r,|t|)^{-1}
\end{equation}
whenever $t\leq 2r$.  Before we turn to this estimate, however, let us first check that
\begin{equation}\label{nn}
|\la_\pm|\cdot \br{|\la_\pm| + |\tau|}^{-1} \leq Cr \cdot \br{|t|+r}^{-1} \quad\quad\text{whenever \,$t-2r\leq \tau\leq
t$.}
\end{equation}

\Case{1} If $|t|\geq 2r$ and $t\geq 0$, then
\begin{equation*}
|\la_\pm| + |\tau| \geq t-\tau -r + |\tau| \geq t-r \geq \frac{t+r}{3} = \frac{|t|+r}{3} \:.
\end{equation*}
Since we also have $|\la_\pm|\leq t-\tau +r\leq 3r$, our claim \eqref{nn} follows.

\Case{2} If $|t|\geq 2r$ and $t\leq 0$, we similarly have
\begin{equation*}
|\la_\pm| + |\tau| \geq t-\tau -r + |\tau| = t-r-2\tau \geq -r-t \geq \frac{r-t}{3} = \frac{r+|t|}{3}
\end{equation*}
and also $|\la_\pm|\leq t-\tau +r\leq 3r$, hence \eqref{nn} follows for this case as well.

\Case{3} If $|t|\leq 2r$ and $r\leq 1$, our claim \eqref{nn} holds because
\begin{equation*}
|\la_\pm|\cdot \br{|\la_\pm| + |\tau|}^{-1} \leq |\la_\pm| \leq 3r \leq 12 r\cdot \br{|t|+r}^{-1}.
\end{equation*}

\Case{4} If $|t|\leq 2r$ and $r\geq 1$, on the other hand, it holds because
\begin{equation*}
|\la_\pm|\cdot \br{|\la_\pm| + |\tau|}^{-1} \leq 1 \leq 4r \cdot \br{|t|+r}^{-1}.
\end{equation*}

Let us now turn to the proof of \eqref{Ipm2}.  Employing \eqref{nn}, we find that
\begin{equation*}
\mathcal{I}_\pm' \leq Cr^a \br{|t|+r}^{-a} \int_{t-2r}^{\min(t,0)} \br{\la_\pm}^{1 -\kappa} \cdot
\br{|\la_\pm|+\tau}^{-\nu} \:d\tau
\end{equation*}
whenever $t\leq 2r$.  In order to establish \eqref{Ipm2}, it thus suffices to show that
\begin{equation}\label{Ipm3}
\mathcal{I}_\pm'' \equiv \int_{t-2r}^t \br{\la_\pm}^{1 -\kappa} \cdot \br{|\la_\pm|+\tau}^{-\nu} \:d\tau \leq
C\br{|t|-r}^{-\nu}.
\end{equation}
Since $\la_+ = t-\tau+r$ is non-negative here, $\mathcal{I}_+''$ is rather easy to treat. To treat $\mathcal{I}_-''$, on
the other hand, it is convenient to write
\begin{align*}
\mathcal{I}_-'' &= \br{t-r}^{-\nu} \int_{t-2r}^{t-r} \br{t-r-\tau}^{1-\kappa} \:d\tau + \int_{t-r}^t \br{\tau +
r-t}^{1-\kappa} \cdot \br{2\tau+r-t}^{-\nu} \:d\tau.
\end{align*}
Using our assumption $\kappa>2$ for the former integral and the substitution $x= 2\tau+r-t$ for the latter, we find
\begin{equation*}
\mathcal{I}_-'' \leq C\br{t-r}^{-\nu} + \frac{1}{2} \,\int_{t-r}^{t+r} \br{\frac{x+r-t}{2}}^{1-\kappa} \cdot
\br{x}^{-\nu} \:dx.
\end{equation*}
Once we now recall that $0<\nu < \nu p<1$ by \eqref{ks}, an application of Lemma \ref{fa5} gives
\begin{equation*}
\mathcal{I}_-'' \leq C\br{t-r}^{-\nu} \leq C\br{|t|-r}^{-\nu}
\end{equation*}
because $|t-r|\geq ||t|-r|$.  This already establishes \eqref{Ipm3}, so the proof is complete.
\end{proof}

\begin{lemma}\label{Jpm}
Let $(r,t)\in \Omega$.  Assuming \eqref{k1}, \eqref{con}, \eqref{ks} and \eqref{iwcs}, one has
\begin{equation*}
\mathcal{J}_\pm \equiv \int_{t-2r}^t |\la_\pm|^{a+1-mp+m} \br{\la_\pm}^{-1} \cdot W_k(|\la_\pm|, |\tau|)^{-p} \:d\tau
\leq Cr^a W_k(r,|t|)^{-1},
\end{equation*}
where $\la_\pm= t-\tau \pm r$, $W_k$ is given by \eqref{W1} and the constant $C$ is independent of $r,t$.
\end{lemma}

\begin{proof}
According to Lemma 3.13 in \cite{Ka2}, we do have the estimate
\begin{equation*}
\int_{\max(t-2r,0)}^t |\la_\pm|^{a+1-mp+m} \br{\la_\pm}^{-1} \cdot W_k(|\la_\pm|, |\tau|)^{-p} \:d\tau \leq Cr^a
W_k(r,|t|)^{-1}
\end{equation*}
whenever $t\geq 0$.  In particular, it suffices to show that
\begin{equation*}
\int_{t-2r}^{\min(t,0)} |\la_\pm|^{a+1-mp+m} \br{\la_\pm}^{-1} \cdot W_k(|\la_\pm|, |\tau|)^{-p} \:d\tau \leq Cr^a
W_k(r,|t|)^{-1}
\end{equation*}
whenever $t\leq 2r$.  Let us now proceed as in the proof of the previous lemma.  Using \eqref{nn}, one may deduce the
last inequality as soon as we know that
\begin{align}\label{Jpm1}
\mathcal{J}_\pm'
&\equiv \int_{t-2r}^t |\la_\pm|^{1-mp+m} \,\br{\la_\pm}^{-1} \cdot \br{|\la_\pm|-\tau}^{-a(p-1)} \cdot \br{|\la_\pm|+
\tau}^{-\nu p} \:d\tau \notag \\
&\leq C\br{|t|-r}^{-\nu}.
\end{align}
In what follows, we only concern ourselves with $\mathcal{J}_-'$, as $\mathcal{J}_+'$ is easier to treat.  Write
\begin{align*}
\mathcal{J}_-'
&= \br{t-r}^{-\nu p} \int_{t-2r}^{t-r} (t-r-\tau)^{1-mp+m} \,\br{t-r-\tau}^{-1} \cdot \br{t-r-2\tau}^{-a(p-1)} \:d\tau \\
&\quad + \br{t-r}^{-a(p-1)} \int_{t-r}^t (r-t+\tau)^{1-mp+m} \,\br{r-t+\tau}^{-1} \cdot \br{r-t+2\tau}^{-\nu p} \:d\tau.
\end{align*}
Changing variables by $x= t-r-2\tau$ for the former integral and by $x= r-t+2\tau$ for the latter, we now arrive at
\begin{align*}
\mathcal{J}_-' &\leq C\br{t-r}^{-\nu p} \int_{r-t}^{3r-t} (x+t-r)^{1-mp+m} \,\br{\frac{x+t-r}{2}}^{-1} \cdot
\br{x}^{-a(p-1)} \:dx \\
&\quad + C\br{t-r}^{-a(p-1)} \int_{t-r}^{t+r} (x+r-t)^{1-mp+m} \,\br{\frac{x+r-t}{2}}^{-1} \cdot \br{x}^{-\nu p} \:dx.
\end{align*}
In view of \eqref{iwcs}, we may then apply Lemma \ref{fa5} with $b= a(p-1)$ and $b= \nu p$ to get
\begin{equation*}
\mathcal{J}_-' \leq C\br{t-r}^{2-mp+m-\nu p-a(p-1)} \leq C\br{t-r}^{-\nu}
\end{equation*}
using \eqref{ide}.  Since this already implies the desired \eqref{Jpm1}, the proof is complete.
\end{proof}

\begin{lemma}\label{fa3}
Let $w\leq 0$ be arbitrary.  Assuming that $a,\nu>0$, one has
\begin{equation*}
B_1 \equiv \int_{-\infty}^w \frac{\br{y}^{-a-\nu} \,dy}{(w-y)^{1-a}} \leq C(a,\nu) \cdot \br{w}^{-\nu}.
\end{equation*}
\end{lemma}

\begin{proof}
Let us change variables by $z= -y$ and write the given integral as
\begin{equation*}
B_1 = \int_{|w|}^{2|w|+1} \frac{\br{z}^{-a-\nu} \,dz}{(w+z)^{1-a}} + \int_{2|w|+1}^\infty \frac{\br{z}^{-a-\nu}
\,dz}{(w+z)^{1-a}} \:.
\end{equation*}
Since $\br{z}$ is equivalent to $\br{w}$ within the former integral and $w+z$ is equivalent to $\br{z}$ within the
latter, we find that
\begin{equation*}
B_1 \leq C\br{w}^{-a-\nu} \int_{|w|}^{2|w|+1} (w+z)^{a-1} \,dz + C \int_{2|w|+1}^\infty \br{z}^{-1-\nu} \,dz.
\end{equation*}
In view of our assumptions that $a,\nu>0$, the desired estimate is now easy to deduce.
\end{proof}

\begin{lemma}\label{fa4}
Let $(r,t)\in \Omega$.  Assuming \eqref{k1}, \eqref{con} and \eqref{ks} with $a\leq 1$, one has
\begin{equation*}
B_2 \equiv \int_{-\infty}^{-|t-r|} \int_0^{\la_-} \frac{\la^m\br{\la}^{1-\kappa} \cdot W_k(\la,|\tau|)^{-1}}{\la_-^m\,
(\la_- -\la)^{1-a}} \:\:d\la\,d\tau \leq C\br{t-r}^{-\nu}
\end{equation*}
and also
\begin{equation}\label{J2n}
B_3 \equiv \int_{-\infty}^{-|t-r|} \int_0^{\la_-} \frac{\la^{1-mp+2m} \cdot W_k(\la,|\tau|)^{-p}}{\la_-^m \,(\la_-
-\la)^{1-a}} \:\:d\la\,d\tau \leq C\br{t-r}^{-\nu},
\end{equation}
where $\la_-= t-\tau- r$, $W_k$ is given by \eqref{W1} and the constant $C$ is independent of $r,t$.
\end{lemma}

\begin{proof}
We shall only prove our second assertion, as the proof of our first assertion is quite similar. Since $\la\leq \la_-$
within the region of integration, we trivially get
\begin{equation*}
B_3 \leq \int_{-\infty}^{-|t-r|} \int_0^{t-r-\tau} \frac{\la^{1-mp+m} \cdot \br{\la-\tau}^{-ap} \br{\la+\tau}^{-\nu
p}}{(t-r-\la -\tau)^{1-a}} \:\:d\la\,d\tau
\end{equation*}
and we may switch to characteristic coordinates $x=\la-\tau$, $y=\la+\tau$ to arrive at
\begin{equation*}
B_3 \leq C\int_{-\infty}^{t-r} \frac{\br{y}^{-\nu p}}{(t-r-y)^{1-a}} \int_{\max(|y|,|t-r|)}^\infty (x+y)^{1-mp+m} \cdot
\br{x}^{-ap} \:dx\,dy.
\end{equation*}
Once we now invoke Lemma \ref{fa2} to treat the inner integral, we find
\begin{align*}
B_3 &\leq C\int_{-\infty}^{t-r} \frac{\br{y}^{-\nu p}}{(t-r-y)^{1-a}} \cdot \max(\br{y}, \br{t-r})^{\nu p -\nu-a} \:dy \\
&= C\int_{-\infty}^{-|t-r|} \frac{\br{y}^{-a-\nu}}{(t-r-y)^{1-a}} \:\:dy + C\br{t-r}^{\nu p-\nu -a} \int_{-|t-r|}^{t-r}
\frac{\br{y}^{-\nu p}}{(t-r-y)^{1-a}} \:\:dy.
\end{align*}
Using Lemma \ref{fa3} for the former integral and Lemma \ref{fa1} for the latter, we get \eqref{J2n}.
\end{proof}

\begin{lemma}\label{I2s}
Let $r,t\in \R$.  Assuming \eqref{ka}, \eqref{con} and \eqref{ks}, one has
\begin{equation*}
\mathbb{I} \equiv \int_{-2|t-r|-1}^{\min(t-r,0)} \int_0^{\la_-} \frac{\la^m \br{\la}^{1-\kappa}
\:d\la\,d\tau}{\la_-^{m+a} \,(\la_- -\la)^{1-a}} \leq C
\end{equation*}
and also
\begin{equation*}
\mathbb{J} \equiv \int_{-2|t-r|-1}^{\min(t-r,0)} \int_0^{\la_-} \frac{\la^{1-mp+2m} \:\:d\la\,d\tau}{\la_-^{m+a} \,(\la_-
-\la)^{1-a}} \leq C\br{t-r}^{2-mp+m},
\end{equation*}
where $\la_-= t-\tau-r$ and the constant $C$ is independent of $r,t$.
\end{lemma}

\begin{proof}
We only concern ourselves with our second assertion, as our first assertion is easier to establish.  Changing variables
by $y= \la_- = t-r-\tau$, let us write
\begin{align*}
\mathbb{J}
&= \int_{\max(t-r,0)}^{1+t-r+2|t-r|} \int_{y/2}^y \frac{\la^{1-mp+2m} \:d\la\,dy}{y^{m+a} \,(y-\la)^{1-a}} +
\int_{\max(t-r,0)}^{1+t-r+2|t-r|} \int_0^{y/2} \frac{\la^{1-mp+2m} \:d\la\,dy}{y^{m+a} \,(y-\la)^{1-a}} \\
&\equiv \mathbb{J}_1 + \mathbb{J}_2.
\end{align*}
When it comes to the former integral, $\la$ is equivalent to $y$ and we easily get
\begin{align*}
\mathbb{J}_1 &\leq C\int_0^{1+3|t-r|} y^{1-mp+m-a} \int_{y/2}^y \,(y-\la)^{a-1} \:d\la\,dy \\
&\leq C\int_0^{1+3|t-r|} y^{1-mp+m} \:dy \leq C\br{t-r}^{2-mp+m}
\end{align*}
because $a>0$ by \eqref{ks} and $2-mp+m>0$ by \eqref{con}.  When it comes to the latter integral, on the other hand,
$y-\la$ is equivalent to $y$, so we find
\begin{equation*}
\mathbb{J}_2 \leq C\int_0^{1+3|t-r|} y^{-m-1} \int_0^{y/2} \la^{1-mp+2m} \:d\la\,dy.
\end{equation*}
Since $2-mp+2m\geq 2-mp+m>0$ by \eqref{ks} and \eqref{con}, the desired estimate follows.
\end{proof}

\begin{lemma}\label{J2}
Let $(r,t)\in \Omega$.  Assuming \eqref{k1}, \eqref{ka}, \eqref{con} and \eqref{ks} with $m\geq 1\geq a$, one has
\begin{equation*}
\mathcal{I}_2 \equiv \int_{-\infty}^{t-r} \int_0^{\la_-} \frac{\la^m \br{\la}^{1-\kappa} \cdot
W_k(\la,|\tau|)^{-1}}{\la_-^m \la_+^a \,(\la_- -\la)^{1-a}} \:\:d\la\,d\tau \leq CW_k(r,|t|)^{-1}
\end{equation*}
and also
\begin{equation*}
\mathcal{J}_2 \equiv \int_{-\infty}^{t-r} \int_0^{\la_-} \frac{\la^{1-mp+2m} \cdot W_k(\la,|\tau|)^{-p}}{\la_-^m \la_+^a
\,(\la_- -\la)^{1-a}} \:\:d\la\,d\tau \leq CW_k(r,|t|)^{-1},
\end{equation*}
where $\la_\pm = t-\tau\pm r$, $W_k$ is given by \eqref{W1} and the constant $C$ is independent of $r, t$.
\end{lemma}

\begin{proof}
We shall only prove our second assertion, as the proof of our first assertion is quite similar.  According to Lemma 3.10
in \cite{Ka2}, we do have the estimate
\begin{equation*}
\int_{\min(t-r,0)}^{t-r} \int_0^{\la_-} \frac{\la^{1-mp+2m} \cdot W_k(\la,|\tau|)^{-p}}{\la_-^m \la_+^a \,(\la_-
-\la)^{1-a}} \:\:d\la\,d\tau \leq CW_k(r,|t|)^{-1}
\end{equation*}
when $t\geq r$.  Since this estimate holds trivially when $t\leq r$, it remains to show that
\begin{equation}\label{sh3}
\mathcal{J}_2' \equiv \int_{-\infty}^{\min(t-r,0)} \int_0^{\la_-} \frac{\la^{1-mp+2m} \cdot W_k(\la,|\tau|)^{-p}}{\la_-^m
\la_+^a \, (\la_- -\la)^{1-a}} \:\:d\la\,d\tau \leq CW_k(r,|t|)^{-1}.
\end{equation}

\Case{1} When $t\geq 0$, it is convenient to express the last integral as the sum of
\begin{equation}\label{J2a}
\mathcal{J}_{21}' \equiv \int_{-\infty}^{-|t-r|} \int_0^{\la_-} \frac{\la^{1-mp+2m} \cdot W_k(\la,|\tau|)^{-p}}{\la_-^m
\la_+^a \,(\la_- -\la)^{1-a}} \:\:d\la\,d\tau
\end{equation}
and
\begin{equation}\label{J2b}
\mathcal{J}_{22}' \equiv \int_{-|t-r|}^{\min(t-r,0)} \int_0^{\la_-} \frac{\la^{1-mp+2m} \cdot
W_k(\la,|\tau|)^{-p}}{\la_-^m \la_+^a \,(\la_- -\la)^{1-a}} \:\:d\la\,d\tau.
\end{equation}

Let us treat $\mathcal{J}_{21}'$ first.  If it happens that $t+r\geq 1$, then $\la_+= t+r-\tau\geq C\br{t+r}$ within the
region of integration.  Using this fact along with Lemma \ref{fa4}, we deduce the desired
\begin{align*}
\mathcal{J}_{21}' &\leq C\br{t+r}^{-a} \int_{-\infty}^{-|t-r|} \int_0^{\la_-} \frac{\la^{1-mp+2m} \cdot
W_k(\la,|\tau|)^{-p}}{\la_-^m \,(\la_- -\la)^{1-a}} \:\:d\la\,d\tau \\
&\leq C\br{t+r}^{-a} \br{t-r}^{-\nu} = CW_k(r,t)^{-1}.
\end{align*}
Assume now that $t+r\leq 1$.  Since we still have $\la_+ = t+r-\tau\geq t+r+1$ for each $\tau\leq -1$, our argument above
ensures that
\begin{align*}
\mathcal{J}_{21}' &\leq CW_k(r,t)^{-1} + \int_{-1}^{-|t-r|} \int_0^{\la_-} \frac{\la^{1-mp+2m} \cdot
W_k(\la,|\tau|)^{-p}}{\la_-^m \la_+^a \,(\la_- -\la)^{1-a}} \:\:d\la\,d\tau.
\end{align*}
Besides, $|t-r|\leq t+r\leq 1$ here, so we need only show that the right hand side is bounded.  Using a trivial estimate
and then Lemma \ref{I2s}, we arrive at the desired
\begin{align*}
\mathcal{J}_{21}' &\leq C + \int_{-2|t-r|-1}^{\min(t-r,0)} \int_0^{\la_-} \frac{\la^{1-mp+2m}
\:\:d\la\,d\tau}{\la_-^{m+a} \,(\la_- -\la)^{1-a}} \leq C.
\end{align*}

Next, we treat $\mathcal{J}_{22}'$.  In treating this integral \eqref{J2b}, we may assume that $t\geq r$ and write
\begin{equation}\label{J2e}
\mathcal{J}_{22}' = \int_{r-t}^0 \int_0^{\la_-} \frac{\la^{1-mp+2m} \cdot W_k(\la,|\tau|)^{-p}}{\la_-^m \la_+^a \,(\la_-
-\la)^{1-a}} \:\:d\la\,d\tau.
\end{equation}
Using a trivial estimate and then Lemma \ref{I2s}, one easily settles the case $t+r\leq 1$ as above.  Assume now that
$t+r\geq 1$.  Since $\la_\pm = t\pm r-\tau \geq t\pm r$ within the region of integration, we find
\begin{equation*}
\mathcal{J}_{22}' \leq \frac{C}{(t-r)^m \br{t+r}^a} \int_{r-t}^0 \int_0^{t-r-\tau} \frac{\la^{1-mp+2m} \cdot
\br{\la-\tau}^{-ap} \br{\la+\tau}^{-\nu p}}{(t-r-\la-\tau)^{1-a}} \:\:d\la\,d\tau.
\end{equation*}
Switching to characteristic coordinates $x=\la-\tau$ and $y= \la+\tau$, we thus find
\begin{equation*}
\mathcal{J}_{22}' \leq \frac{C}{(t-r)^m \br{t+r}^a} \int_{r-t}^{t-r} \frac{\br{y}^{-\nu p}}{(t-r-y)^{1-a}}
\int_{|y|}^{3(t-r)} (x+y)^{1-mp+2m} \,\br{x}^{-ap} \:dx\,dy.
\end{equation*}
Note that $2(x+1) \geq 2x\geq x+y$ here, while \eqref{con} implies
\begin{equation*}
2-mp+2m-ap > m-a \geq 0
\end{equation*}
because we are assuming that $m\geq a$.  Combining these facts, we now get
\begin{equation*}
\mathcal{J}_{22}' \leq \frac{C(t-r)^{2-mp+m-ap}}{\br{t+r}^a} \int_0^{t-r} \frac{\br{y}^{-\nu p}}{(t-r\pm y)^{1-a}}
\:\:dy.
\end{equation*}
Besides, $\nu p<1$ by \eqref{ks}, so we may apply Lemma \ref{fa1} with $z= t-r\geq 0$ to arrive at
\begin{equation*}
\mathcal{J}_{22}' \leq C(t-r)^{2-mp+m-ap+a} \br{t-r}^{-\nu p} \br{t+r}^{-a}.
\end{equation*}
In view of our assumption \eqref{con} and our inequality \eqref{ide}, we also have
\begin{equation*}
0< 2-(m+a)(p-1) \leq \nu (p-1),
\end{equation*}
so we may combine the last two equations to deduce the desired estimate \eqref{sh3}.

\Case{2} When $t\leq 0$, we have $t<r$ and we shall express the integral \eqref{sh3} as the sum of
\begin{equation}\label{J2c}
\mathcal{J}_{21}'' \equiv \int_{-\infty}^{2(t-r)-1} \int_0^{\la_-} \frac{\la^{1-mp+2m} \cdot
W_k(\la,|\tau|)^{-p}}{\la_-^m \la_+^a \,(\la_- -\la)^{1-a}} \:\:d\la\,d\tau
\end{equation}
and
\begin{equation}\label{J2d}
\mathcal{J}_{22}'' \equiv \int_{2(t-r)-1}^{t-r} \int_0^{\la_-} \frac{\la^{1-mp+2m} \cdot W_k(\la,|\tau|)^{-p}}{\la_-^m
\la_+^a \,(\la_- -\la)^{1-a}} \:\:d\la\,d\tau.
\end{equation}

To treat the first part \eqref{J2c}, we note that $\la_+ = t+r -\tau \geq 3r-t+1> r-t+1$ within the region of
integration. In view of Lemma \ref{fa4}, this gives
\begin{equation*}
\mathcal{J}_{21}'' \leq \br{r-t}^{-a} \int_{-\infty}^{t-r} \int_0^{\la_-} \frac{\la^{1-mp+2m} \cdot
W_k(\la,|\tau|)^{-p}}{\la_-^m \,(\la_- -\la)^{1-a}} \:\:d\la\,d\tau \leq C\br{r-t}^{-a-\nu},
\end{equation*}
which implies the desired \eqref{sh3} because $|r+t|\leq r+|t| = r-t$ whenever $t\leq 0$.

To treat the second part \eqref{J2d}, we use the fact that $\la_+\geq \la_-$ to trivially obtain
\begin{equation*}
\mathcal{J}_{22}'' \leq \int_{2(t-r)-1}^{t-r} \int_0^{\la_-} \frac{\la^{1-mp+2m} \cdot \br{\la-\tau}^{-ap} \br{\la
+\tau}^{-\nu p}}{\la_-^{m+a} \,(\la_- -\la)^{1-a}} \:\:d\la\,d\tau.
\end{equation*}
Since $\la-\tau\geq -\tau\geq r-t$ and $\la+\tau\leq t-r< 0$ within the region of integration, we get
\begin{equation*}
\mathcal{J}_{22}'' \leq \br{r-t}^{-ap-\nu p} \int_{2(t-r)-1}^{t-r} \int_0^{\la_-} \frac{\la^{1-mp+2m}
\:\:d\la\,d\tau}{\la_-^{m+a} \,(\la_- -\la)^{1-a}} \leq C\br{r-t}^{2-mp+m-ap-\nu p}
\end{equation*}
by Lemma \ref{I2s}.  In view of our inequality \eqref{ide}, the desired estimate \eqref{sh3} follows.
\end{proof}

\section{Existence of solutions}\label{ex}
In this section, we prove the two main results of this paper. Before we do that, however, let us first combine the
estimates of the previous section and give the

\begin{proof_of}{Theorem \ref{bes}}
There are two estimates we need to verify and we shall only verify the first, as the proof of the second is similar.  To
estimate the Duhamel integral \eqref{mb1}, we apply Proposition \ref{Les}.  Given a function $G\in \mathcal{C}^1(\Omega)$
that satisfies the singularity condition
\begin{equation}\label{sg3}
G(\la,\tau) = O(\la^{-2m-2+\de}) \quad\quad \text{as\, $\la\to 0$}
\end{equation}
for some fixed $\de>0$, Proposition \ref{Les} provides the estimate
\begin{align}\label{hg}
|D^\b [\mathscr{L}G](r,t)|
&\leq Cr^{j-|\b|-m-a} \int_{-\infty}^t \int_{|\la_-|}^{\la_+} \frac{\la^{m-j+1}}{(\la-\la_-)^{1-a}} \cdot \sum_{s=0}^j
\la^s \,|\d_\la^s G(\la,\tau)| \:d\la\,d\tau \notag\\
&\quad + Cr^{j-|\b|-m} \int_{-\infty}^{t-r} \int_0^{\la_-} \frac{\la^{2m-j+1}}{\la_-^m \la_+^a \,(\la_- -\la)^{1-a}}
\cdot \sum_{s=0}^j \la^s \,|\d_\la^s G(\la,\tau)| \:d\la\,d\tau \notag\\
&\quad + Cr^{j-|\b|-m-a} \int_{t-2r}^t  |\la_\pm|^{a+m-j+1} \left[ \sum_{s=0}^{j-1} \la^s \,|\d_\la^s
G(\la,\tau)| \right]_{\la= |\la_\pm|} d\tau
\end{align}
whenever $(r,t)\in \Omega$ and $|\b|\leq j\leq 1$.

We now use this fact with $G(\la,\tau)= F(u(\la,\tau))$.  By Lemma \ref{bas} with $j=j_0=0$, we have
\begin{equation*}
G(\la,\tau)= F(u(\la,\tau)) = O(\la^{-mp}) = O(\la^{-2m-2+\de_1}) \quad\quad \text{as\:\:$\la\to 0$,}
\end{equation*}
where $\de_1= 2-mp+2m$ is positive by \eqref{con}.  In particular, our estimate \eqref{hg} does hold for the special case
$G= F(u)$. Besides, the sums that appear in the right hand side are those of Lemma \ref{bas}, according to which
\begin{equation*}
\sum_{s=0}^{j_0} \la^s \,|\d_\la^s F(u(\la,\tau))| \leq 2Ap ||u||^p \cdot\la^{j-mp} \br{\la}^{j_0-j} \cdot
W_k(\la,|\tau|)^{-p}, \quad\quad j_0= j, j-1.
\end{equation*}
Once we now insert this fact in \eqref{hg}, we obtain an estimate of the form
\begin{equation*}
|D^\b [\mathscr{L} F(u)](r,t)| \leq C||u||^p \cdot r^{j-|\b|-m} \cdot (r^{-a} \mathcal{J}_1 + \mathcal{J}_2 + r^{-a}
\mathcal{J}_\pm),
\end{equation*}
where $\mathcal{J}_1$, $\mathcal{J}_2$ and $\mathcal{J}_\pm$ are as in Lemmas \ref{J1}, \ref{J2} and \ref{Jpm},
respectively.  The assumptions we imposed in those lemmas are not different from the ones imposed in this theorem, except
for the inequality $a\leq 1\leq m$ that appears in Lemma \ref{J2}.  Since our definition \eqref{am} shows that $a\leq
1\leq m$ whenever $n\geq 4$, however, we may employ Lemmas \ref{J1}, \ref{Jpm} and \ref{J2} to get
\begin{equation}\label{an1}
|D^\b [\mathscr{L} F(u)](r,t)| \leq C||u||^p \cdot r^{j-|\b|-m} \cdot W_k(r,|t|)^{-1}, \quad\quad |\b|\leq j\leq 1.
\end{equation}
Finally, we claim that this also implies our first assertion \eqref{mb1}, namely that
\begin{equation*}
|D^\b [\mathscr{L} F(u)](r,t)| \leq C||u||^p \cdot r^{1-|\b|-m} \br{r}^{|\b|-1} \cdot W_k(r,|t|)^{-1}, \quad\quad
|\b|\leq 1.
\end{equation*}
Indeed, if $r\leq 1$, one may obtain the last inequality through the special case $j=1$ of \eqref{an1}.  If $r\geq 1$, on
the other hand, one may obtain it through the special case $j=|\b|\leq 1$.
\end{proof_of}

\begin{proof_of}{Theorem \ref{et}}
Our iteration argument is quite standard and also similar to the one used in \cite{Ka2, Kb1}, so we only give a sketch of
the proof.  As we have already noted, one may decrease the decay rates $k,\kappa$ to ensure that \eqref{k1} through
\eqref{ka} hold.  By \eqref{kl} and \eqref{k2}, our decay rate $k$ is then such that
\begin{equation*}
\frac{n-1}{2} < k < \frac{n+1}{2} \leq n-1
\end{equation*}
because $n\geq 3$.  This allows us to invoke Lemma \ref{dec} to obtain a unique solution $u_0^-\in X$ to the homogeneous
problem \eqref{he} such that $||u_0^-||\leq C_0\e$.  Setting $u_0= u_0^-$, we then recursively define a sequence
$\{u_i\}$ using the formula
\begin{equation}\label{it}
u_{i+1} = u_0^- + \mathscr{L}F(u_i) - \mathscr{L}(Vu_i), \quad\quad i\geq 0.
\end{equation}
In view of Theorem \ref{bes}, we do have the estimate
\begin{align*}
|\d_r^j \mathscr{L}F(u_i)| &\leq C_1||u_i||^p \cdot r^{1-j-m} \br{r}^{j-1} \cdot W_k(r,|t|)^{-1}, \quad\quad j=0,1
\end{align*}
whenever $u_i\in X$.  In view of the definition \eqref{N} of our norm, we thus have
\begin{equation*}
||\mathscr{L}F(u_i)|| \leq C_1 ||u_i||^p.
\end{equation*}
Since we may similarly estimate $||\mathscr{L}(Vu_i)||$, equation \eqref{it} leads us to
\begin{equation*}
||u_{i+1}|| \leq ||u_0^-|| + C_1||u_i||^p + C_1V_0||u_i|| \leq C_0\e + C_1||u_i||^p + C_1V_0||u_i||.
\end{equation*}
As long as $\e$ and $V_0$ are sufficiently small, we may now use the approach of \cite{Ka2, Kb1} to deduce that our
sequence converges to an element $u\in X$ with small norm.  In view of \eqref{it}, this gives us a solution to the
integral equation
\begin{equation}\label{ie}
u - u_0^- = \mathscr{L}F(u) - \mathscr{L}(Vu).
\end{equation}
Using another application of Theorem \ref{bes}, we thus arrive at
\begin{align*}
|D^\b (u-u_0^-)| &\leq C_1||u||^p \cdot r^{1-|\b|-m} \br{r}^{|\b|-1} \cdot W_k(r,|t|)^{-1} \\
&\quad + C_1V_0 ||u|| \cdot r^{1-|\b|-m} \br{r}^{|\b|-1} \cdot W_k(r,|t|)^{-1},
\end{align*}
where $D= (\d_r,\d_t)$ and $|\b|\leq 1$.  This trivially implies the first part of the theorem.

To prove the second part of the theorem, we define
\begin{equation*}
u_0^+(r,t) = u_0^-(r,t) + \int_{-\infty}^\infty [LG(u(\cdot\,,\tau))](r,t-\tau) \:d\tau,
\end{equation*}
where $L$ is the Riemann operator and we have set $G(u) = F(u) - V(r)\cdot u$.  Then $u_0^+$ is clearly a solution to the
homogeneous equation \eqref{he}.  Recalling the integral equation \eqref{ie} and the definition \eqref{Du} of the Duhamel
operator, we can rewrite the last equation as
\begin{equation*}
u_0^+(r,t)= u(r,t) + \int_t^\infty [LG(u(\cdot\,,\tau))](r,t-\tau) \:d\tau.
\end{equation*}
Since the Riemann operator $L$ is odd in $t$, a change of variables then easily leads us to
\begin{equation*}
u(r,t) - u_0^+(r,t) = \int_{-\infty}^{-t} \:[LG(u(\cdot\,,-\sigma))](r,-t-\sigma) \:d\sigma.
\end{equation*}
Besides, the right hand side bears a close resemblance to the Duhamel operator \eqref{Du}, so one can estimate this
expression in exactly the same way that we used to estimate \eqref{ie}.
\end{proof_of}

\begin{proof_of}{Corollary \ref{sc}}
The two estimates that we need to prove are quite similar to one another, so let us merely focus on the first one, namely
\begin{equation}\label{fe}
||u(\cdot\,,t) - u_0^-(\cdot\,,t)||_e \leq C||u||\cdot \br{t}^{-\theta}, \quad\quad t\leq 0.
\end{equation}
As it is well-known, this estimate will follow immediately once we know that
\begin{align*}
\int_{-\infty}^t \left( \int_0^\infty F(u(r,\tau))^2 \cdot r^{n-1}\:dr\right)^{1/2} d\tau
&\leq C||u||\cdot \br{t}^{-\theta}, \quad\quad t\leq 0; \\
\int_{-\infty}^t \left( \int_0^\infty V(r)^2 \cdot u(r,\tau)^2 \cdot r^{n-1}\:dr\right)^{1/2} d\tau
&\leq C||u||\cdot \br{t}^{-\theta}, \quad\quad t\leq 0.
\end{align*}
Needless to say, these two estimates are also similar to one another.  Let us merely establish the former by showing that
\begin{equation}\label{dn1}
H(\tau) \equiv \int_0^\infty F(u(r,\tau))^2 \cdot r^{n-1}\:dr \leq C||u||^2\cdot \br{\tau}^{-2\theta-2}, \quad\quad
\tau\leq 0
\end{equation}
for some $\theta>0$.  Using our assumption \eqref{F} on $F$ together with the fact that $u$ belongs to the Banach space
\eqref{X}, we get
\begin{equation*}
|F(u(r,\tau))| \leq A |u(r,\tau)|^p \leq A||u||^p \cdot r^{p(1-m)} \br{r}^{-p} \cdot \br{|\tau|+r}^{-ap}
\br{|\tau|-r}^{-\nu p}.
\end{equation*}
Next, we insert this fact in \eqref{dn1}.  Since $2(a+m)= n-1$ by \eqref{am}, we find that
\begin{equation*}
H(\tau) \leq C||u||^{2p} \int_0^\infty r^{2a+2m+2p(1-m)} \br{r}^{-2p} \cdot \br{|\tau|+r}^{-2ap} \br{|\tau|-r}^{-2\nu p}
\:dr.
\end{equation*}
Since $||u||^{2p} \leq ||u||^2$ whenever $||u||\leq 1$, it suffices to show that
\begin{equation}\label{dn15}
\widetilde{H}(\tau) \equiv \int_0^\infty r^{2a+2m+2p(1-m)} \br{r}^{-2p} \cdot \br{|\tau|+r}^{-2ap} \br{|\tau|-r}^{-2\nu
p} \:dr \leq C\br{\tau}^{-2\theta-2}
\end{equation}
whenever $\tau\leq 0$.  Now, if $r\leq 1$, then $\br{|\tau|\pm r}$ is equivalent to $\br{\tau}$, so we easily get
\begin{align*}
\widetilde{H}(\tau) &\leq C\br{\tau}^{-2p(a+\nu)} \int_0^1 r^{2a+2m-2mp+2p} \:dr \\
&\quad + C\int_0^\infty \br{r}^{2a+2m-2mp} \cdot \br{|\tau|+r}^{-2ap} \br{|\tau|-r}^{-2\nu p} \:dr.
\end{align*}
To see that the former integral is finite, we note that
\begin{equation*}
2a+2m-2mp+2p > 2ap + 2p-4 \geq 3p-4>-1
\end{equation*}
by \eqref{con} and since $a=1/2,1$.  This allows us to arrive at
\begin{align}
\widetilde{H}(\tau)
&\leq C\br{\tau}^{-2p(a+\nu)} + C\br{\tau}^{-2p(a+\nu)} \int_0^{|\tau|/2} \br{r}^{2a+2m-2mp} \:dr \notag\\
&\quad +C\int_{2|\tau|}^\infty \br{r}^{2a+2m-2p(a+m+\nu)} \:dr
+ C\br{\tau}^{-2(a+m)(p-1)} \int_{|\tau|/2}^{2|\tau|} \br{|\tau|-r}^{-2\nu p} \:dr \notag\\
&\equiv \widetilde{H}_1(\tau) + \widetilde{H}_2(\tau) + \widetilde{H}_3(\tau) + \widetilde{H}_4(\tau) \label{dn2}.
\end{align}

Focusing on the first two terms for the moment, it is clear that we have
\begin{equation*}
\widetilde{H}_1(\tau) + \widetilde{H}_2(\tau) \leq C\br{\tau}^{-2p(a+\nu)} +C_\de \br{\tau}^{2a+2m-2p(a+m+\nu)+1+\de}
\end{equation*}
for an arbitrary choice of $\de>0$.  Since $a+\nu = k-m$ by \eqref{ks}, we thus have
\begin{align*}
\widetilde{H}_1(\tau) + \widetilde{H}_2(\tau)
&\leq C\br{\tau}^{-2p(k-m)} +C_\de \br{\tau}^{-2\nu +2k-2kp+1+\de} \\
&\leq C\br{\tau}^{-2p(k-m)} +C_\de \br{\tau}^{-2\nu -3+\de}
\end{align*}
because $k(p-1)\geq 2$ by \eqref{k1}.  Arguing similarly for the third term in \eqref{dn2}, we find
\begin{equation*}
\widetilde{H}_3(\tau) = C\int_{2|\tau|}^\infty \br{r}^{-2\nu+2k-2kp} \:dr \leq C\int_{2|\tau|}^\infty \br{r}^{-2\nu-4}
\:dr \leq C\br{\tau}^{-2\nu-3}
\end{equation*}
because $\nu>0$ by \eqref{ks}.  As for the fourth term in \eqref{dn2}, we get
\begin{align*}
\widetilde{H}_4(\tau)
&\leq C\br{\tau}^{-2(a+m)(p-1)} +C_\de \br{\tau}^{2a+2m-2p(a+m+\nu)+1+\de} \\
&\leq C\br{\tau}^{-2(a+m)(p-1)} +C_\de \br{\tau}^{-2\nu -3+\de}
\end{align*}
exactly as before.  Once we now combine these observations, we obtain the estimate
\begin{equation*}
\widetilde{H}(\tau)\leq C\br{\tau}^{-2p(k-m)} + C\br{\tau}^{-2(a+m)(p-1)} + C_\de \br{\tau}^{-2\nu -3+\de}
\end{equation*}
for an arbitrary choice of $\de>0$.  Taking $\de= 3-2a$, we finally deduce the desired
\begin{equation*}
\widetilde{H}(\tau)\leq C\br{\tau}^{-2p(k-m)} + C\br{\tau}^{-2(a+m)(p-1)} + C_\de \br{\tau}^{-2(k-m)} \leq
C\br{\tau}^{-2\theta-2}
\end{equation*}
with $\theta$ defined as in \eqref{th}.  Since $\theta>0$ whenever $k>m+1$, the desired \eqref{dn15} follows.
\end{proof_of}

\end{document}